\documentclass[english]{article}
\usepackage{amsfonts}
\usepackage{amsmath}
\usepackage{geometry}
\begin{document}

\title{Bessel and Struve Related Integrals}
\author{Bernard J. Laurenzi \\
Department of Chemistry, University at Albany, The State University of New York \\
1400 Washington Ave., Albany N. Y. 12222}
\date{January 6, 2014}
\maketitle

\begin{abstract}
Analytic expressions for integrals which arise in a theory of atomic
structure due to Schwinger and Englert are evaluated in terms of Bessel and
Struve functions.

\end{abstract}

\bigskip

Integrals which arise within a theory of atomic structure due to Schwinger
and Englert \cite{SandE} are investigated in the work below. \ In keeping
with the spirit of their approach, analytical expressions for those integral
are presented below. \ We denote the integrals in question by%
\[
\mathcal{S}(z,\zeta )=\int_{0}^{\pi /2}\cos \theta \,\sin ^{2}\theta \,\sin
(z\cos \theta )\,\sin (\varsigma \cos ^{2}\theta )\,d\theta , 
\]%
\[
\mathcal{C}(z,\zeta )=\int_{0}^{\pi /2}\cos \theta \,\sin ^{2}\theta \,\cos
(z\cos \theta )\cos (\varsigma \cos ^{2}\theta )\,d\theta . 
\]%
These integrals will be shown to involve higher derivatives of the Bessel
functions of the first kind i.e. $J_{\nu }(z)$ and of the Struve functions $%
\mathbf{H}_{\nu }\left( z\right) $ \cite{NBS} respectively. The derivatives
in turn can be reduced to expressions containing the corresponding Bessel
and Struve functions of orders zero and one.

\bigskip

\section{The $\mathcal{S}(z,\protect\zeta )$ Integral}

In the case of the ``\textit{sine''} integral $\mathcal{S}(z,\zeta )$ we
begin by expanding the term $\sin (\varsigma \cos ^{2}\theta )$ in a power
series in $\varsigma $ to get 
\begin{equation}
\mathcal{S}(z,\zeta )=\sum_{\kappa =0}^{\infty }\frac{(-1)^{\kappa
}\varsigma ^{2\,\kappa +1}}{(2\kappa +1)!}\int_{0}^{\pi /2}\cos ^{4\,\kappa
+3}\theta \,\sin ^{2}\theta \,\sin (z\cos \theta )\,d\theta .  \label{eq1}
\end{equation}%
From the theory of Bessel functions \cite{Soni} we have the related integral%
\[
J_{1}(z)/z=\frac{2}{\pi }\int_{0}^{\pi /2}\cos (z\cos \theta )\sin
^{2}\theta \,d\theta . 
\]%
Differentiation\thinspace of this relation with respect to $z,$ $4\kappa +3$
times produces the integrals occurring in (1). \ We have%
\[
\mathcal{S}(z,\zeta )=\frac{\pi }{2}\sum_{\kappa =0}^{\infty }\frac{%
(-1)^{\kappa }\varsigma ^{2\,\kappa +1}}{(2\kappa +1)!}\frac{d\,^{4\,\kappa
+3}\left[ J_{1}(z)/z\right] }{dz^{\,4\,\kappa +3}}, 
\]%
an expression which requires analytic expressions for the higher derivatives
of $J_{1}(z)/z.$

\bigskip

\subsection{The Higher Derivatives of $J_{1}(z)/z$}

Once again, using an identity from the theory of Bessel functions of the
first kind \cite{NBS2} i.e.%
\[
\left( \frac{1}{z}\frac{d}{dz}\right) ^{k}\frac{J_{\nu }(z)}{z^{\nu }}%
=(-1)^{k}\frac{J_{\,\nu +k}(z)}{z^{\,\nu +k}}, 
\]%
we will make use of that relation which in this application is more usefully
rewritten as%
\begin{equation}
\frac{d^{\,k}}{d\,(z^{2\,})^{k}}\left( \frac{J_{1}(z)}{z}\right) =\left( -%
\frac{1}{2}\right) ^{k}\frac{J_{\,k+1}(z)}{z^{\,k+1}}.  \label{eq2}
\end{equation}%
In this form equation (2) makes it possible to obtain the required higher
derivatives. \ Using Fa\`{a} di Bruno's \cite{Faa} formula for the
differentiation of a function of a function (i.e. a composite function) we
have 
\[
\frac{d^{\,k}}{d\,z^{\,}{}^{k}}\left( \frac{J_{1}(z)}{z^{{}}}\right)
=k!\sum_{i=0}^{{\Large \lfloor }\frac{k}{2}{\Large \rfloor }}\frac{%
(2z)^{k-2i}}{i!\,(k-2i)!}\frac{d^{\,k-i}[J_{1}(z)/z]}{d\,\left(
z^{2}\right) ^{\,}{}^{k-i}}. 
\]%
Using this together with equation (2), the sought after expression for the
higher derivatives of $J_{1}(z)/z$ is 
\begin{equation}
\frac{d^{\,k}}{d\,z^{\,}{}^{k}}\left( \frac{J_{1}(z)}{z^{{}}}\right) =%
{\small 2(-1)}^{k}{\small \,k!}\sum_{i=0}^{{\Large \lfloor }\frac{k}{2}%
{\Large \rfloor }}\frac{(-1)^{i}}{i!\,\,(k-2i)!}\frac{J_{k+1-i}(z)}{%
(2z)^{i+1}},  \label{eq3}
\end{equation}%
where $\lfloor z\rfloor $ i.e. the \textit{floor} of $z,$ is the largest
integer $\leq $ $z.$ \ The Bessel functions of order $k+1-i$ on the
right-hand side of equation (3) are reducible to expressions which contains
the products of each of the functions i.e. $J_{0}(z),$ and $J_{1}(z)$ with
corresponding polynomials containing powers of $1/z.$\bigskip

Using recurrence relations (for decreasing distant neighbors) for the Bessel
functions \cite{Brych} i.e.%
\[
J_{\nu }\left( z\right) =C_{n}\left( \nu ,z\right) \,J_{\nu -n}\left(
z\right) -C_{n-1}(\nu ,z)\,J_{\nu -n-1}\left( z\right) , 
\]%
with%
\[
C_{0}\left( \nu ,z\right) =1,C_{1}\left( \nu ,z\right) =2(\nu -1)/z, 
\]%
\[
C_{n}\left( \nu ,z\right) =\frac{2(\nu -n)}{z}C_{n-1}\left( \nu ,z\right)
-C_{n-2}\left( \nu ,z\right) , 
\]%
we have for this application the required \textit{reduced} recurrence
relation 
\begin{equation}
J_{\nu }\left( z\right) =C_{\nu -1}\left( \nu ,z\right) \,J_{1}\left(
z\right) -C_{\nu -2}(\nu ,z)\,J_{0}\left( z\right) .  \label{eq4}
\end{equation}%
The difference equation for $C_{n}\left( \nu ,z\right) $ has solutions%
\[
C_{n}(\nu ,z)=(-2/z)^{n}(1-\nu )_{n}\,\cdot \,_{2}F_{3}\tbinom{\frac{1-n}{2}%
,\;-\frac{n}{2};\;-z^{2}}{1-\nu ,\;-n,\;\nu -n}, 
\]%
where $(a)_{n\text{ }}$is the Pochhammer symbol and $_{2}F_{3}$ is a
generalized hypergeometric function. \ \bigskip

In the case of the $J_{0}(z)$ prefactor i.e. $C_{\nu -2}(\nu ,z),$ the
hypergeometric function is reducible in the case of integer $\nu $ (after
some algebra), to a polynomial in $1/z$ i.e.%
\[
C_{\nu -2}(\nu ,z)=\sum_{j=0}^{\lceil \frac{\nu }{2}-1\rceil }\frac{\Gamma
(\nu -j)\,\Gamma (\nu -1-j)}{\Gamma (\nu -1-2j)}\frac{(-1)^{\,j\,+1}}{%
j!\,(j+1)!}\left( \frac{2}{z}\right) ^{\nu -2-2j}=\mathcal{R}^{(0)}(\nu
,\,z), 
\]%
where $\lceil z\rceil $ is the ceiling of $z$, i.e. the smallest integer $%
\geq $ to $z.$ The polynomials $C_{\nu -2}(\nu ,z)$ are in fact special
cases of the well known Lommel polynomials i.e. $R_{\nu -2,2}(z)$ \cite%
{Watson} . \ For simplicity in notation we will denote these special cases
of the Lommel polynomials by $\mathcal{R}^{(0)}(\nu ,\,z)$.

\bigskip

In a similar way, the prefactor of $J_{1}(z)$ i.e. $C_{\nu -1}(\nu ,z)$
reduces to 
\[
C_{\nu -1}(\nu ,z)=\sum_{j=0}^{\lceil \frac{\nu -1}{2}\rceil }\frac{\Gamma
(\nu -j)^{2}}{\Gamma (\nu -2j)}\frac{(-1)^{j}}{j!^{2}}\left( \frac{2}{z}%
\right) ^{\nu -1-2j}=\mathcal{R}^{\left( 1\right) }(\nu ,z), 
\]%
which are the Lommel polynomials $R_{\nu -1,1}(z)$ and will be referred to
as $\mathcal{R}^{\left( 1\right) }(\nu ,\,z).$ Using (3,4) the $k^{th}$%
derivative of $J_{1}(z)/z$ is then given by 
\[
\frac{d^{\,k}}{d\,z^{\,}{}^{k}}\left( \frac{J_{1}(z)}{z^{{}}}\right)
=(-1)^{k}\,\left[ \mathcal{P}_{1}(k,z)\,J_{1}(z)\,-\mathcal{P}%
_{0}(k,z)\,J_{0}(z)\right] \,, 
\]%
where $\mathcal{P}_{1}(k,z)$ and $\mathcal{P}_{0}(k,z)$ are polynomials in $%
1/z$ i.e. \ 
\[
\mathcal{P}_{1}(k,z)=2k!
\sum_{i=0}^{{\Large \lfloor }\frac{k}{2}{\Large \rfloor }}\frac{(-1)^{i}}{i!\,(k-2i)!}\frac{\mathcal{R}^{\left( 1\right)}(k+1-i\,,\,z)}{(2\,z)^{i+1}}, 
\]%
\[
\mathcal{P}_{0}(k,z)=2k!\sum_{i=0}^{{\Large \lfloor }\frac{k}{2}{\Large %
\rfloor }}\frac{(-1)^{i}}{i!\,(k-2i)!}\frac{\mathcal{R}^{\left( 0\right)
}(k+1-i\,,\,z)}{(2\,z)^{i+1}}. 
\]

The $J_{1}(z)/z$ derivatives are seen to be oscillatory functions of $z$
with small and rapidly decreasing amplitudes. \ Direct calculation of the $%
J_{1}(z)/z$ derivatives shows that they either vanish at small $z$ in the
case of odd values of $k$ or for even $k=2m$ have values $(-1)^{m}\Gamma
(m+1/2){\large /}2\sqrt{\pi }(m+1)!$. \ \bigskip

The integral $\mathcal{S}(z,\zeta )$ is then given by%
\[
\mathcal{S}(z,\zeta )=\frac{\pi }{2}\sum_{\kappa =0}^{\infty }\frac{%
(-1)^{\kappa }\varsigma ^{2\kappa +1}}{(2\kappa +1)!}\left[ \mathcal{P}%
_{0}(4\kappa +3,z)\,J_{0}(z)-\mathcal{P}_{1}(4\kappa +3,z)\,J_{1}(z)\right]
, 
\]%
where%
\begin{eqnarray*}
\mathcal{P}_{0}(4\kappa +3,z) &=&2(4\kappa +3)!\sum_{i=0}^{{\Large 2\kappa
+1}}\frac{(-1)^{i}}{i!\,(4\kappa +3-2i)!}\frac{\mathcal{R}^{\left(
0\right) }(4\kappa +4-i\,,\,z)}{(2\,z)^{i+1}}, \\
\mathcal{P}_{1}(4\kappa +3,z) &=&2(4\kappa +3)!\sum_{i=0}^{{\Large 2\kappa
+1}}\frac{(-1)^{i}}{i!\,(4\kappa +3-2i)!}\frac{\mathcal{R}^{\left(
1\right) }(4\kappa +4-i\,,\,z)}{(2\,z)^{i+1}}.
\end{eqnarray*}

\bigskip

\section{The $\mathcal{C}(z,\protect\zeta )$ Integral}

The ``\textit{cosine''} integral $\mathcal{C}(z,\zeta )$ is treated in a
manner similar to the ``sine'' integral $\mathcal{S}(z,\zeta )$ discussed
above i.e.%
\begin{equation}
\mathcal{C}(z,\zeta )=\sum_{\kappa =0}^{\infty }\frac{(-1)^{\kappa
}\varsigma ^{2\kappa }}{(2\kappa )!}\int_{0}^{\pi /2}\cos ^{4\kappa
+1}\theta \,\sin ^{2}\theta \,\cos (z\cos \theta )\,d\theta .  \label{eq5}
\end{equation}%
Here however, the analysis is considerably more complicated due to the
occurrence of the Struve functions \cite{Struve}. The integral
representation for the Struve function $\mathbf{H}_{1}(z)/z$ 
\[
\mathbf{H}_{1}(z)/z=\frac{2\,}{\pi }\int_{0}^{\pi /2}\sin (z\cos \theta
)\sin ^{2}\theta \,d\theta , 
\]%
upon differentiation $4\kappa +1$ times with respect to $z$ gives the
integrals occurring in (5) with the result that 
\[
\mathcal{C}(z,\zeta )=\frac{\pi }{2}\sum_{\kappa =0}^{\infty }\frac{%
(-1)^{\kappa }\varsigma ^{2\,\kappa }}{(2\kappa )!}\frac{d\,^{4\,\kappa +1}%
\left[ \mathbf{H}_{1}(z)/z\right] }{dz^{4\,\kappa +1}}. 
\]

\subsection{The Higher Derivatives of $\mathbf{H}_{1}(z)/z\hspace{4in}$}

In the case of the Struve functions, relations analogous to those in
equation (2) do not exist. \ However, using the recurrence relation \cite%
{Struve}%
\[
\mathbf{H}_{1}(z)=\frac{2}{\pi }-\mathbf{H}_{-1}(z), 
\]%
we get%
\begin{equation}
\frac{d^{\,k}}{d\,z^{\,}{}^{k}}\left( \frac{\mathbf{H}_{1}(z)}{z^{{}}}%
\right) =\frac{2k!}{\pi }\frac{(-1)^{k}\,}{z^{k+1}}-\frac{d^{\,k}}{%
d\,z^{\,k}}\left( \frac{\mathbf{H}_{-1}(z)}{z^{{}}}\right) .  \label{eq6}
\end{equation}%
The identity \cite{NIST} 
\[
\left( \frac{1}{z^{{}}}\frac{d}{dz}\right) ^{k}\frac{\mathbf{H}_{-1}(z)}{%
z^{{}}}=\frac{\mathbf{H}_{\,-k-1}(z)}{z^{\,k+1}}, 
\]%
analogous to the one involving the Bessel function used above will be useful
in this case. \ Rewritten as%
\[
\frac{d^{\,\,k}}{d\,(z^{2}){}^{k}}\left( \frac{\mathbf{H}_{-1}(z)}{z^{{}}}%
\right) =\left( \frac{1}{2}\right) ^{k}\frac{\mathbf{H}_{-k\,-1}(z)}{z^{k+1}}%
, 
\]%
it will prove helpful here. \ Using the Fa\`{a} di Bruno formula as employed
above, together with the latter relation we get 
\[
\frac{d^{\,k}}{d\,z^{\,}{}^{k}}\left( \frac{\mathbf{H}_{-1}(z)}{z^{{}}}%
\right) =2\,k!\sum_{i=0}^{{\Large \lfloor }\frac{k}{2}{\Large \rfloor }}%
\frac{1}{i!\,(k-2i)!}\frac{\mathbf{H}_{-\,k-1+i}(z)}{(2\,z)^{i+1}}. 
\]%
Using (6), we have the required derivatives i.e. 
\begin{equation}
\frac{d^{\,k}}{d\,z^{\,}{}^{k}}\left( \frac{\mathbf{H}_{1}(z)}{z^{{}}}%
\right) =\frac{2\,k!}{\pi ^{{}}}\frac{(-1)^{k}\,}{z^{k+1}}-\,{\small 2}\,k%
{\small !}\sum_{i=0}^{{\Large \lfloor }\frac{k}{2}{\Large \rfloor }}\frac{1%
}{i!\,(k-2i)!}\frac{\mathbf{H}_{-\,k-1+i}(z)}{(2\,z)^{i+1}}.  \label{eq7}
\end{equation}%
The Struve functions occurring in (7) can also be reduced to expressions
containing only $\mathbf{H}_{0}(z)$ and $\mathbf{H}_{1}(z).$ To do that, we
note as a first step that the Struve functions with negative integer orders
are related to those with the corresponding positive orders by the relations %
\cite{appnote}%
\[
\mathbf{H}_{-\,\nu }(z)=(-1)^{\nu }\,\mathbf{H}_{\nu }(z)+\sum_{j=0}^{\nu -1}%
\frac{(-1)^{\,\,j}}{\Gamma (\,j+3/2)\,\Gamma (\,j+3/2-\nu )}\left( \frac{z%
}{2}\right) ^{2\,j+1-\nu }. 
\]%
As a result (7) becomes%
\begin{eqnarray}
\frac{d^{\,k}}{d\,z^{\,}{}^{k}}\left( \frac{\mathbf{H}_{1}(z)}{z^{{}}}%
\right) &=&{\small 2\,k!\,(-1)}^{k}\left\{ \frac{1}{\pi \,\,z^{\,k+1}}%
+\sum_{i=0}^{\lfloor k/2\rfloor }\frac{\left( -1\right) \,^{i}}{i!\,(k-2i)!}%
\frac{\mathbf{H}_{\,k+1-i}(z)}{(2\,z)^{\,i+1}}\right\}  \label{eq8} \\
&&-{\small k!}\left( \frac{2}{z}\right) ^{k+1}\sum_{i=0}^{\lfloor
k/2\rfloor }\frac{\left( 1/2\right) \,^{2\,i+1}}{i!\,(k-2i)!}%
\sum_{j=0}^{k-i}\frac{(-1)^{\,j}\,\left( z/2\right) ^{2\,j}}{\Gamma
(j+3/2)\,\Gamma (j+1/2+i-k)}.  \nonumber
\end{eqnarray}%
The last (double) sum in equation (8) can be rewritten in ascending powers
of $z/2$ as%
\[
-k!\left( \frac{2}{z}\right) ^{k+1}{\LARGE [}\sum_{j=0}^{k-\,\lfloor
k/2\rfloor -1}\frac{(-1)^{\,\,j}\,\left( z/2\right) ^{2\,j}}{\Gamma
(j+3/2)\,}\sum_{i=0}^{\lfloor k/2\rfloor }\frac{\left( 1/2\right)
\,^{2\,i+1}}{i!\,(k-2i)!\Gamma (j+1/2+i-k)} 
\]%
\[
+\sum_{j=\,k-\,\lfloor k/2\rfloor }^{k}\frac{(-1)^{\,\,j}\,\left(
z/2\right) ^{2\,j}}{\Gamma (j+3/2)\,}\sum_{i=0}^{j-k}\frac{\left(
1/2\right) \,^{2\,i+1}}{i!\,(k-2i)!\Gamma (j+1/2+i-k)}{\LARGE ].} 
\]%
The first of these two sum makes a contribution only when $j=0,$ all other
terms vanishing. \ The surviving term from that sum cancels the leading term
in (8) and we get%
\[
\frac{d^{\,k}}{d\,z^{\,}{}^{k}}\left( \frac{\mathbf{H}_{1}(z)}{z^{{}}}%
\right) = 
\]%
\begin{equation}
{\small \,k!\,(-1)}^{k}\left[ 2\sum_{i=0}^{\lfloor k/2\rfloor }\frac{\left(
-1\right) \,^{i}}{i!\,(k-2i)!}\frac{\mathbf{H}_{\,k+1-i}(z)}{(2\,z)^{\,i+1}}%
-\sum_{j=0}^{\lfloor k/2\rfloor }\frac{(-1)^{\,\,j}\,\left( z/2\right)
^{k-1-2\,j}}{\Gamma (k+3/2-j)\,}\sum_{i=0}^{j}\frac{\left( 1/2\right)
\,^{2\,i+1}}{i!\,(k-2i)!\Gamma (i+1/2-j)}\right] .  \label{eq9}
\end{equation}

\bigskip

As in the case of the ``sine'' integral $\mathcal{S}(z,\zeta )$ it is
necessary to use recurrence relations (for decreasing distant neighbors) for
the Struve functions $\mathbf{H}_{\,k+1-i}(z)$ in order to rewrite those
functions in terms of $\mathbf{H}_{0}(z)$ and\thinspace $\mathbf{H}_{1}(z).$
The Wolfram Function Site a product of \textit{Mathematica} \cite{Brych2}
has given the expressions (due to Yu A. Brychkov) needed for the reduction
of the $\mathbf{H}_{\,k+1-i}(z)$ functions i.e. 
\begin{equation}
\mathbf{H}_{\nu }(z)=C_{\nu -1}(\nu ,z)\,\mathbf{H}_{1}(z)-C_{\nu -2}(\nu
,z)\,\mathbf{H}_{0}(z)+\mathbf{S(}\nu ,z\mathbf{),}  \label{eq10}
\end{equation}%
where%
\[
\mathbf{S(}\nu ,z\mathbf{)=}\frac{1}{\sqrt{\pi }}\sum_{j=0}^{\nu -2}\frac{%
(z/2)^{\nu -1-j}}{\Gamma (\nu +1/2-j)}C_{j}(\nu ,z), 
\]%
\begin{equation}
C_{j}(\nu ,z)=(-2/z)^{\,j}(1-\nu )_{j}\,\cdot \,_{2}F_{3}\tbinom{\frac{1-j}{2%
},\;-\frac{j}{2};\;-z^{2}}{1-\nu ,\;-j,\;\nu -j}.  \label{eq11}
\end{equation}%
This is a similar but more complicated set of relations than those
encountered in the case of the $J_{\nu }(z)$ functions. \ In spite of that,
we see that the prefactors of $\mathbf{H}_{0}(z)$ and $\mathbf{H}_{1}(z)$ in
(10) are identical to the corresponding ones for $J_{0}(z)$ and $J_{1}(z)$
in the reduced recurrence relations for $J_{\nu }(z).$ In this case however,
an additional term i.e. $\mathbf{S(}\nu ,z\mathbf{)}$, occurs and requires
special attention.

In the case of the $\mathbf{H}_{0}(z)$ prefactor i.e. $C_{\nu -2}(\nu ,z),$
the corresponding hypergeometric function (as was seen above in the case of
the Besssel functions) reduces \ to a polynomial in $1/z$ i.e.%
\[
C_{\nu -2}(\nu ,z)=\sum_{\kappa =0}^{\lceil \frac{\nu }{2}-1\rceil }\frac{%
\Gamma (\nu -\kappa )\,\Gamma (\nu -1-\kappa )}{\Gamma (\nu -1-2\kappa )}%
\frac{(-1)^{\kappa +1}}{\kappa !\,(\kappa +1)\,!}\left( \frac{2}{z}\right)
^{\nu -2-2\,\kappa }=R_{\nu -2,\,2}(z). 
\]%
In a similar way the prefactor of $\mathbf{H}_{1}(z)$ i.e. $C_{\nu -1}(\nu
,z)$ reduces to 
\[
C_{\nu -1}(\nu ,z)=\sum_{\kappa =0}^{\lceil \frac{\nu -1}{2}\rceil }\frac{%
\Gamma (\nu -\kappa )^{2}}{\Gamma (\nu -2\kappa )}\frac{(-1)^{\kappa }}{%
\kappa !^{2}}\left( \frac{2}{z}\right) ^{\nu -1-2\,\kappa }=R_{\nu
-1,\,1}(z). 
\]%
A table containing the first eight of the Lommel polynomials $R_{\nu
-1,1}(z),$ $R_{\nu -2,2}(z)$ have been given below in the appendix.

\subsection{The sum $\mathbf{S(}\protect\nu ,z\mathbf{)=}\frac{1}{\protect%
\sqrt{\protect\pi }}\sum_{j=0}^{\protect\nu -2}\frac{(z/2)^{\protect\nu %
-1-j}}{\Gamma (\protect\nu +1/2-j)}C_{j}(\protect\nu ,z)\hspace{5in}$}

For either even or odd values of $\nu $ it is possible to simplify the term $%
\mathbf{S(}\nu ,z\mathbf{)}$ in (10). \ Substitution of the expression for $%
C_{j}(\nu ,z)$ in (11) into $\mathbf{S(}\nu ,z\mathbf{)}$ and simplifying
the result gives 
\[
\mathbf{S(}\nu ,z\mathbf{)}=\frac{1}{\pi }\sum_{j=0}^{\nu -2}\frac{%
\,\Gamma (\nu )\,}{\Gamma (2\nu -2\,j)}2^{2\nu -1-2j}\,(\frac{z}{2})^{\nu
-1-j}\,_{2}F_{3}\tbinom{\frac{1-j}{2},\;-\frac{j}{2};\;-z^{2}}{1-\nu
\;,\;-j,\;\nu -j}. 
\]

\subsection{The subcases $\mathbf{S(}\protect\nu ,z\mathbf{)}$ with $\protect%
\nu =2m$, $\protect\nu =2m+1$ for integer $m\hspace{4in}$}

In the case where $\nu =2m$ we have%
\[
\mathbf{S(}2m,z\mathbf{)}=\frac{1}{\pi }\sum_{j=0}^{2m-2}\frac{\,\Gamma
(2m)\,}{\Gamma (4m-2\,j)}2^{4m-1-2j}\,(\frac{z}{2})^{2m-1-j}\,_{2}F_{3}%
\tbinom{\frac{1-j}{2},\;-\frac{j}{2};\;-z^{2}}{1-2m,\;-j,\;2m-j}. 
\]%
Separating the sum over $j$ into its even and odd terms we have \ 
\begin{eqnarray}
\mathbf{S(}2m,z\mathbf{)} &=&\frac{1}{\pi }\sum_{i=0}^{m-1}\frac{\,\Gamma
(2m)\,}{\Gamma (4m-4i)}2^{4m-1-4i}\,(\frac{z}{2})^{2m-1-4i}\,_{2}F_{3}%
\tbinom{-i+\frac{1}{2},\,-i\,;\,-z^{2}}{1-2m,\,-2i,\,2m-2i}  \label{eq12} \\
&&+\frac{1}{\pi }\sum_{i=0}^{m-2}\frac{\,\Gamma (2m)\,}{\Gamma (4m-4i-2)}%
2^{4m-3-4i}\,(\frac{z}{2}{\small )}^{2m-3-4i}\,_{2}F_{3}\tbinom{-i,\,-i-%
\frac{1}{2};\,-z^{2}}{1-2m,\,-1-2i,\,2m-2i-1}.  \nonumber
\end{eqnarray}%
The hypergeometric functions in the two sum in equation (12) reduce to 
\[
{\small \,}_{2}{\small F}_{3}\tbinom{-i+\frac{1}{2},\,-i\,;\,-z^{2}}{%
1-2m,\,-2i,\,2m-2i}=\frac{\Gamma (2m-2i)}{\Gamma (2m)}\sum_{\kappa =0}^{i}%
\frac{(-1)^{\kappa }\,\Gamma (2m-\kappa )\,\Gamma (2i+1-\kappa )}{\Gamma
(2i+1-2\kappa )\,\Gamma (\kappa +2m-2i)}\frac{(z/2)^{2\kappa }}{\kappa !}, 
\]%
and 
\[
\,\,_{2}{\small F}_{3}\tbinom{-i,\,-i-\frac{1}{2};\,-z^{2}}{%
1-2m,\,-1-2i,\,2m-2i-1}=\frac{\Gamma (2m-1-2i)}{\Gamma (2m)}\sum_{\kappa
=0}^{i}\frac{(-1)^{\kappa }\,\Gamma (2m-\kappa )\,\Gamma (2i+2-\kappa )}{%
\Gamma (2i+2-2\kappa )\,\Gamma (\kappa +2m-1-2i)}\frac{(z/2)^{2\kappa }}{%
\kappa !}, 
\]%
respectively. \ Recombining these terms we get for $\mathbf{S(}2m,z\mathbf{)}
$ the expression%
\[
\mathbf{S(}2m,z\mathbf{)}=\frac{2^{4m-1}}{\pi }\sum_{j=0}^{2m-2}\frac{%
\Gamma (2m-2j)\,2^{-2j}}{\Gamma (4m-2j)}\sum_{\kappa =0}^{\lceil \frac{j-1}{%
2}\rceil }\frac{(-1)^{\kappa }\,\Gamma (2m-\kappa )\,\Gamma (j+1-\kappa )}{%
\Gamma (j+1-2\kappa )\,\Gamma (\kappa +2m-j)}\frac{(z/2)^{2\kappa -2j+2m-1}%
}{\kappa !} 
\]

A similar analysis of the sum $\mathbf{S(}2m+1,z\mathbf{)}$ allows one to
conclude that in each case, a single relation valid for arbitrary $j$ and $%
\nu $ is obtainable i.e. 
\[
_{2}F_{3}\tbinom{\frac{1-j}{2},\,\,-\frac{j}{2};\,\,-z^{2}}{1-\nu
,\,\,-j,\,\,\nu -j}=\frac{\Gamma (\nu -j)}{\Gamma (\nu )}\sum_{\kappa
=0}^{\lceil \frac{j\,-1}{2}\rceil }\frac{(-1)^{\kappa }\,\Gamma (\nu
-\kappa )\,\Gamma (j+1-\kappa )}{\Gamma (j+1-2\kappa )\,\Gamma (\kappa +\nu
-j)}\frac{(z/2)^{2\kappa }}{\kappa !}. 
\]%
As a result, 
\begin{equation}
\mathbf{S(}\nu ,z\mathbf{)}=\frac{2^{2\,\nu -1}}{\pi }\sum_{j=0}^{\nu -2}%
\frac{\,\Gamma (\nu -j)\,\,2^{-2j}}{\Gamma (2\nu -2\,j)}\,\sum_{\kappa
=0}^{\lceil \frac{j\,-1}{2}\rceil }\frac{(-1)^{\kappa }\,\Gamma (\nu
-\kappa )\,\Gamma (j+1-\kappa )}{\Gamma (j+1-2\,\kappa )\,\Gamma (\kappa
+\nu -j)}\frac{(z/2)^{2\,\kappa -2\,j+\,\nu -1}}{\kappa !}.  \label{eq13}
\end{equation}%
If the order of summation in the double sum in equation (13) is
interchanged, that expression (for $\nu \geq 2$) can be written in ascending
powers of $2/z$ as 
\begin{eqnarray*}
\mathbf{S(}\nu ,z\mathbf{)} &=&\frac{1}{\pi }(\frac{z}{2})^{\nu
-1}\sum_{\mu =0}^{\nu -2-\,{\large \lceil }\frac{\nu \,-1}{2}{\large \rceil }%
}\frac{\Gamma (\mu +1/2)\,}{\Gamma (\nu +1/2-\mu )}(\frac{2}{z})^{2\mu } \\
&&+\frac{1}{\sqrt{\pi }}(\frac{z}{2})^{\nu -1}\sum_{\mu =\,\nu -1-\,%
{\large \lceil }\frac{\nu \,-1}{2}{\large \rceil }}^{\nu -2}\frac{\mu !\,}{%
\Gamma (\nu -\mu )}\left[ \sum_{i=0}^{\nu -\mu -2}\frac{(-1)^{\,i}\,\Gamma
(\nu -i)}{\Gamma (\nu +1/2-\mu -i)\,\Gamma (\mu +1-i)\,\,i!}\right] (\frac{2%
}{z})^{2\mu }.
\end{eqnarray*}%
This can be rewritten as 
\begin{eqnarray}
\mathbf{S(}\nu ,z\mathbf{)} &=&\frac{1}{\pi }(\frac{z}{2})^{\nu
-1}\sum_{\mu =0}^{{\large \lfloor }\frac{\nu -1}{2}{\large \rfloor \,-1}}%
\frac{\Gamma (\mu +1/2)\,}{\Gamma (\nu +1/2-\mu )}(\frac{2}{z})^{2\mu }
\label{eq14} \\
&&+\frac{1}{\sqrt{\pi }}(\frac{z}{2})^{3-\nu }\sum_{\mu =\,0}^{{\large %
\lceil }\frac{\nu -1}{2}{\large \rceil \,-1}}\left( \sum_{i=0}^{\mu }\frac{%
(-1)^{\,i}\,\Gamma (\nu -i)}{\Gamma (\nu +5/2-i)\,\Gamma (\nu -1-\mu
-i)\,\,i!}\right) \frac{(\nu -2-\mu )!}{\Gamma (\mu +2)}(\frac{z}{2}%
)^{2\mu },  \nonumber
\end{eqnarray}%
where we have used the relation 
\[
\lceil N/2\rceil +\lfloor N/2\rfloor =N, 
\]%
valid for integers (which follows from the properties \cite{ceil} of the
ceiling and floor functions) in the first summation's limits and we have
reindexed the second i.e. the double sum over $\mu $. \ Gathering the terms
found above, the final expression for $\mathbf{H}_{\nu }(z)$\thinspace is 
\begin{eqnarray}
\mathbf{H}_{\nu }(z) &=&\,\mathbf{H}_{1}(z)\sum_{\kappa =0}^{\lceil \frac{%
\nu \,-1}{2}\rceil }\frac{\Gamma (\nu -\kappa )^{2}}{\Gamma (\nu -2\,\kappa
)}\frac{(-1)^{\kappa }}{\kappa !^{2}}\left( \frac{2}{z}\right) ^{\nu
-1-2\,\kappa }  \label{eq15} \\
&&-\mathbf{H}_{0}(z)\sum_{\kappa =0}^{\lceil \frac{\nu }{2}-1\rceil }\frac{%
\Gamma (\nu -\kappa )\,\Gamma (\nu -1-\kappa )}{\Gamma (\nu -1-2\kappa )}%
\frac{(-1)^{\kappa +1}}{\kappa !\,(\kappa +1)!}\left( \frac{2}{z}\right)
^{\nu -2-2\,\kappa }+\mathbf{S(}\nu ,z\mathbf{)},  \nonumber
\end{eqnarray}%
with $\mathbf{S(}\nu ,z\mathbf{)}$ given in (14). \bigskip

Finally, using (9, 15) the expression for the higher derivative of\ $\mathbf{%
H}_{1}(z)/z$ is given by

\begin{eqnarray}
\frac{(-1)^{k}}{k!}\frac{d^{\,k}}{d\,z^{\,}{}^{k}}\left( \frac{\mathbf{H}%
_{1}(z)}{z^{{}}}\right) &=&\left\{ \sum_{i=0}^{\lceil \frac{k}{2}\rceil }%
\frac{(-1/4)^{i}}{i!\,(k-2\,i)!}\frac{1}{2}\sum_{\kappa =0}^{\lceil \frac{%
k\,-\,i}{2}\rceil }\frac{\Gamma (k+1-i-\kappa )^{2}}{\Gamma
(k+1-i-2\,\kappa )}\frac{(-1)^{\,\kappa }}{\left( \kappa !\right) ^{2}}%
\left( \frac{z}{2}\right) ^{2\,\kappa -k-1}\right\} \mathbf{H}_{1}(z)
\label{eq16} \\
&&-\left\{ \sum_{i=0}^{\lceil \frac{k}{2}\rceil }\frac{(-1/4)^{i}}{%
i!\,(k-2i)!}\frac{1}{2}\sum_{\kappa =0}^{\lceil \frac{k-1-i}{2}\rceil }%
\frac{\Gamma (k+1-i-\kappa )\Gamma (k-i-\kappa )}{\,\Gamma (k-i-2\kappa )}%
\frac{(-1)^{\,\kappa }}{\kappa !(\kappa +1)!}\left( \frac{z}{2}\right)
^{2\,\kappa -k}\right\} \mathbf{H}_{0}(z)  \nonumber \\
&&+2\sum_{i=0}^{\lceil \frac{k}{2}\rceil }\frac{(-1)^{i}}{i!\,(k-2i)!}%
\frac{\mathbf{S(}k+1-i,\,z\mathbf{)}}{(2z)^{i+1}}-\sum_{j=0}^{\lfloor 
\frac{k}{2}\rfloor }\frac{(-1)^{\,j}(z/2)^{k-1-2\,j}}{\Gamma (k+3/2-j)}%
\sum_{i=0}^{j}\frac{\left( \frac{1}{2}\right) ^{2\,i+1}}{i!(k-2i)!\,\Gamma
(i+1/2-j)}.  \nonumber
\end{eqnarray}%
The prefactors of $\mathbf{H}_{1}(z)$ and $\mathbf{H}_{0}(z)$ in equation
(16) are now seen to be double sums, and can be simplified when the orders
of the summations are interchanged. \ Furthermore, when $\mathbf{S(}k+1-i,z%
\mathbf{)}$ as represented in (14) is substituted into (16) and the orders
of summation are interchanged, internal cancellation of many of the
resulting terms occurs and we obtain%
\[
{\small (-1)}^{k}\frac{d^{\,\,k}}{d\,z{}^{k}}\left( \frac{\mathbf{H}_{1}(z)%
}{z^{{}}}\right) =\mathbf{H}_{0}(z)\;{\LARGE \sigma }_{0}(k,z)\cdot
(2/z)^{k}+\mathbf{H}_{1}(z)\;{\LARGE \sigma }_{1}(k,z)\cdot (2/z)^{k+1}+%
{\LARGE \sigma }_{2}(k,z)\cdot (2/z)^{k-1}, 
\]%
where ${\LARGE \sigma }_{0}(k,z),{\LARGE \sigma }_{1}(k,z),$and ${\LARGE %
\sigma }_{2}(k,z)$ are polynomials in $z$\ i.e.%
\begin{eqnarray*}
{\LARGE \sigma }_{0}(k,z) &=&k!{\Huge [}\sum_{\nu =0}^{\lceil \frac{1}{2}%
\lfloor \frac{k}{2}\rfloor \,-\frac{1}{2}\,\rceil }\left( \frac{1}{2}%
\sum_{i=0}^{\lceil \frac{k}{2}\rceil }\frac{({\small -\,}1/4)^{i}}{%
i!\,(k-2i)!}\frac{(k-\nu -i)!\,(k-1-\nu -i)!}{(k-1-2\,\nu -i)!}\right) 
\frac{(-1)^{\nu +1}}{\nu !\,(\nu +1)!}\left( \frac{z}{2}\right) ^{2\nu } \\
&&+\sum_{\nu =\,\lceil \frac{1}{2}\lfloor \frac{k}{2}\rfloor \,+\frac{1}{2}%
\,\rceil }^{\lceil \frac{k-1}{2}\,\rceil }\left( \frac{1}{2}%
\sum_{i=0}^{k-1-2\nu }\frac{({\small -\,}1/4)^{i}}{i!\,(k-2i)!}\frac{%
(k-\nu -i)!\,(k-1-\nu -i)!}{(k-1-2\,\nu -i)!}\right) \frac{(-1)^{\nu +1}}{%
\nu !\,(\nu +1)!}\left( \frac{z}{2}\right) ^{2\nu }{\Huge ]},
\end{eqnarray*}%
\begin{eqnarray*}
{\LARGE \sigma }_{1}(k,z) &=&k!{\Huge [}\sum_{\nu =0}^{\lceil \frac{1}{2}%
\lfloor \frac{k}{2}\rfloor \,\rceil }\left( \frac{1}{2}\sum_{i=0}^{\lceil 
\frac{k}{2}\rceil }\frac{({\small -\,}1/4)^{i}}{i!(k-2i)!}\frac{(k-\nu
-i)!^{2}}{(k-2\nu -i)!}\right) \frac{(-1)^{\nu }}{\left( \nu !\right) ^{2}}%
\left( \frac{z}{2}\right) ^{2\nu } \\
&&+\sum_{\nu =1+\lceil \frac{1}{2}\lfloor \frac{k}{2}\rfloor \,\rceil
}^{\lceil \frac{k}{2}\,\rceil }\left( \frac{1}{2}\sum_{i=0}^{k+1-2\nu }%
\frac{({\small -\,}1/4)^{i}}{i!(k-2i)!}\frac{(k-\nu -i)!^{2}}{(k-2\nu -i)!}%
\right) \frac{(-1)^{\nu }}{\left( \nu !\right) ^{2}}\left( \frac{z}{2}%
\right) ^{2\nu }{\Huge ]},
\end{eqnarray*}

\begin{eqnarray*}
{\LARGE \sigma }_{2}(k,z) &=&k!{\Huge [}\frac{({\small -\,}\frac{1}{4}%
)^{\lfloor \frac{k}{2}\rfloor +1}\delta \lbrack 1,\,\lceil \frac{k}{2}\rceil
\,-\,\lfloor \frac{k}{2}\rfloor \,]}{\Gamma (\,\lfloor \frac{k}{2}\rfloor \,+%
\frac{3}{2})\,\Gamma (\,\lceil \frac{k}{2}\rceil \,+\frac{3}{2})}(\frac{z}{2%
})^{k-1}+\frac{1}{2\pi }\sum_{\nu =1}^{\lceil \frac{1}{2}\lfloor \frac{k}{2%
}\rfloor \,\rceil }\frac{\mathfrak{C}(k,\,\,\lceil \frac{k}{2}\rceil ,\,\nu
)}{\nu !}\left( \frac{z}{2}\right) ^{2\nu -2} \\
&&+\frac{1}{2\pi }\sum_{\nu =1+\lceil \frac{1}{2}\lfloor \frac{k}{2}%
\rfloor \,\rceil }^{\lceil \frac{k}{2}\,\rceil }\left\{ \frac{\mathfrak{C}%
(k,\,\,k+1-2\nu ,\,\nu )}{\nu !}+\sum_{i=k+2-2\nu }^{\lceil \frac{k}{2}%
\,\rceil }\frac{({\small -}1/4)^{i}}{i!\,(k-2i)!}\frac{\Gamma (k+1/2-\nu
-i)}{\Gamma (\nu +3/2)}\right\} \left( \frac{z}{2}\right) ^{2\nu -2}{\Huge ]%
},
\end{eqnarray*}%
where $\delta \lbrack i,j]$ is the Kronecker delta (introduced to prevent
over counting in a sum) and when it occurs $1/0!$ is counted as zero. \ The
coefficients $\mathfrak{C}(k,a,\nu )$ are defined as 
\[
\mathfrak{C}(k,a,\nu )=\sum_{j=0}^{a}\frac{({\small -}1/4)\,^{j}}{j\,!}%
\frac{(k-\nu -j)!}{(k-2\,j)!}\sum_{i=0}^{\nu -1}\frac{(-1)^{i}}{i\,!}%
\frac{\sqrt{\pi }\,\;(k-j-i)\,!}{(k-\nu -j-i)\,!\,\,\Gamma (\nu +3/2-i)\,}. 
\]

Direct calculation of the derivatives $d^{\,\,k}[\mathbf{H}%
_{1}(z)/z]/d\,z{}^{k}$ shows that these quantities like those involving $%
J_{1}(z)/z$ also behave as oscillatory functions of $z$ with small and
rapidly decreasing amplitudes. \ At $z=0$ the amplitude is zero for even $k$
and for $k=2m+1$ the amplitude is $(-1)^{m}m!{\LARGE /}2\sqrt{\pi }\,\Gamma
(m+5/2).$ \ In the case of both higher derivatives, the finite values of
their amplitudes at the origin are of order $1/m^{3/2}$ for large $m.$\ \
Further, as anticipated in the work above, the prefactors of $\mathbf{H}%
_{0}(z)$ and $\mathbf{H}_{1}(z)$ and those of $J_{0}(z)$ and $J_{1}(z)$ are
seen to be identical. As a result, alternative representations of the
quantities ${\LARGE \sigma }_{i}(k,z)/z^{k+i}$ and $\mathcal{P}_{i}(k,z)$
i.e. as explicit polynomials in $1/z$ or as sums of Lommel polynomials have
been given.

Finally the expression for $\mathcal{C}(z,\zeta )$ is (with some additional
cancellation occurring because of the odd values of $k$ encountered in this
application) 
\[
\mathcal{C}(z,\zeta )=\frac{\pi }{2}\sum_{\kappa =0}^{\infty }\frac{%
(-1)^{\kappa +1}\varsigma ^{2\kappa }}{(2\,\kappa )!}\left[ 
\begin{array}{c}
\mathbf{H}_{0}(z)\;{\LARGE \sigma }_{0}({\small 4\kappa +1},z)\cdot {\small %
(2/z)}^{4\kappa +1}+\mathbf{H}_{1}(z)\;{\LARGE \sigma }_{1}({\small 4\kappa
+1},z)\cdot ({\small 2/z})^{4\kappa +2} \\ 
+\;{\LARGE \sigma }_{2}({\small 4\kappa +1},z)\cdot ({\small 2/z})^{4\kappa }%
\end{array}%
\right] , 
\]%
with 
\begin{eqnarray*}
{\LARGE \sigma }_{0}({\small 4\kappa +1},z) &=&({\small 4\kappa +1)!}{\Huge [%
}\sum_{\nu =0}^{\kappa }\left( \frac{1}{2}\sum_{i=0}^{2\kappa +1}\frac{(%
{\small -\,}1/4)^{i}}{i!\,(4\kappa +1-2i)!}\frac{(4\kappa +1-\nu
-i)!\,(4\kappa -\nu -i)!}{(4\kappa -2\nu -i)!}\right) \frac{(-1)^{\nu +1}}{%
\nu \,!\,(\nu +1)!}\left( \frac{z}{2}\right) ^{2\nu } \\
&&+\sum_{\nu =\kappa +1}^{2\kappa }\left( \frac{1}{2}\sum_{i=0}^{4\kappa
+1-2\nu }\frac{({\small -\,}1/4)^{i}}{i!\,(4\kappa +1-2i)!}\frac{(4\kappa
+1-\nu -i)!\,(4\kappa -\nu -i)!}{(4\kappa -2\nu -i)!}\right) \frac{%
(-1)^{\nu +1}}{\nu \,!\,(\nu +1)!}\left( \frac{z}{2}\right) ^{2\nu }{\Huge ]%
},
\end{eqnarray*}%
\begin{eqnarray*}
{\LARGE \sigma }_{1}({\small 4\kappa +1},z) &=&({\small 4\kappa +1)!}{\Huge [%
}\sum_{\nu =0}^{\kappa }\left( \frac{1}{2}\sum_{i=0}^{{\small 2\kappa +1}}%
\frac{({\small -\,}1/4)^{i}}{i!({\small 4\kappa +1}-2i)!}\frac{({\small %
4\kappa +1}-\nu -i)!^{2}}{({\small 4\kappa +1}-2\nu -i)!}\right) \frac{%
(-1)^{\nu }}{\left( \nu !\right) ^{2}}\left( \frac{z}{2}\right) ^{2\nu } \\
&&+\sum_{\nu =\kappa +1}^{2\kappa }\left( \frac{1}{2}\sum_{i=0}^{{\small %
4\kappa +2}-2\nu }\frac{({\small -\,}1/4)^{i}}{i!(4\kappa +1-2i)!}\frac{(%
{\small 4\kappa +1}-\nu -i)!^{2}}{({\small 4\kappa +1}-2\nu -i)!}\right) 
\frac{(-1)^{\nu }}{\left( \nu !\right) ^{2}}\left( \frac{z}{2}\right)
^{2\nu }{\Huge ]},
\end{eqnarray*}%
\begin{eqnarray*}
{\LARGE \sigma }_{2}({\small 4\kappa +1},z) &=&\frac{({\small 4\kappa +1})!%
}{2\pi }{\Huge [}\sum_{\nu =1}^{\kappa }\frac{\mathfrak{C}({\small 4\kappa
+1},\,\,{\small 2\kappa +1},\,\nu )}{\nu !}\left( \frac{z}{2}\right) ^{2\nu
-2} \\
&&+\sum_{\nu =\,\kappa +1}^{2\kappa }\left\{ \frac{\mathfrak{C}({\small %
4\kappa +1},\,\,{\small 4\kappa +2}-2\nu ,\,\nu )}{\nu !}+\sum_{i=\,{\small %
4\,\kappa +3}-2\nu }^{2\kappa }\frac{({\small -}1/4)^{i}}{i!({\small %
4\kappa +1}-2i)!}\frac{\Gamma ({\small 4\kappa }+3/2-\nu -i)}{\Gamma (\nu
+3/2)}\right\} \left( \frac{z}{2}\right) ^{2\nu -2}{\Huge ]}.
\end{eqnarray*}

\bigskip 

\begin{center}
\appendix Appendix
\end{center}

In the table below we list a few of the Lommel polynomials referred to in
the text.

\[
\begin{tabular}{|c|c|c|}
\hline
$\nu $ & $R_{\nu -1,1}(z)$ & $R_{\nu -2,2}(z)$ \\ \hline
$1$ & $1$ & $z$ \\ \hline
$2$ & $2/z$ & $1$ \\ \hline
$3$ & $8/z^{2}-1$ & $4/z$ \\ \hline
$4$ & $48/z^{3}-8/z$ & $24/z^{2}-1$ \\ \hline
$5$ & $384/z^{4}-72/z^{2}+1$ & $192/z^{3}-12/z$ \\ \hline
$6$ & $3840/z^{5}-768/z^{3}+18/z$ & $1920/z^{4}-144/z^{2}+1$ \\ \hline
$7$ & $46080/z^{6}-9600/z^{4}+288/z^{2}-1$ & $23040/z^{5}-1920/z^{3}+24/z$
\\ \hline
$8$ & $645120/z^{7}-138240/z^{5}+4800/z^{3}-32/z$ & $%
322560/z^{6}-28800/z^{4}+480/z^{2}-1$ \\ \hline
\end{tabular}%
\]

\end{document}